\documentclass[11pt]{article}
\usepackage{graphicx}
\usepackage{amssymb}
\usepackage{amsmath}

\newcommand{\refeq}[1]{\textrm{(\ref{eq:#1})}}    
\newcommand{\refsec}[1]{Section~\ref{sec:#1}}     

\newcommand{\tr}{\mathrm{trace}} 			

\title{Distribution-free factor analysis --- Estimation theory and applicability to high-dimensional data.}
\author{Rolf Sundberg, Stockholm University \\ Uwe Feldmann, University of Saarland}
\begin{document}

\maketitle
\begin{abstract} We here provide a distribution-free approach to the random factor analysis model. We show that it leads to the same estimating equations as for the classical ML estimates under normality, but more easily derived, and valid also in the case of more variables than observations ($p>n$). For this case we also advocate a simple iteration method. In an illustration with $p=2000$ and $n=22$ it was seen to lead to convergence after just a few iterations. We show that there is no reason to expect Heywood cases to appear, and that the factor scores will typically be precisely estimated/predicted as soon as $p$ is large. We state as a general conjecture that the nice behaviour is not despite $p>n$, but because $p>n$.
\end{abstract}

\vspace{1cm}

\vspace{1.5cm}
\section{Introduction}
In this paper we consider parameter estimation in a distribution-free version of the standard (Gaussian) factor analysis (FA) model, with special emphasis on the case of more variables than observations.
The FA model means describing a sample  $x_1,\ldots,x_n$ of $p$-dimensional vectors as 
\begin{equation} \label{eq:Modelforx}
x_i = \mu + \Lambda f_i + e_i \,, \hspace{5mm} i=1,\ldots,n.
\end{equation}
Here $\mu$ is the mean value vector, $\Lambda$ is a $p\times k$ coefficients (loadings) matrix, $k<\min(n,p)$, and the $f_i$s are mutually independent latent $k$-vectors (factor scores), standardized to zero mean and unit covariance matrix $I_k$ (for identifiability). The $e_i$s are assumed mutually independent $p$-vectors with uncorrelated components and diagonal covariance matrix  $\Psi^2$. Also, $f_i$ and $e_i$ should be mutually independent. In matrix form we write \refeq{Modelforx}  as
$X = \mu {\bf1} + F  \Lambda^T + E$, with the vectors of  \refeq{Modelforx} as rows. 

Usually, normality of $f$ and $e$ in  \refeq{Modelforx}  is assumed, and more observations than variables, that is $n>p$. Then Gaussian maximum likelihood methods can be used, and are more or less standard. However, in recent years interest has increased both in more robust methods and in methods for the case of more variables than observations, $p>n$. Among papers having appeared after the comprehensive review by Bartholomew \& Knott (1999, ch. 3), we mention Robertson \& Symons (2007), who study extension of Gaussian maximum likelihood to the case $p>n$, and a number of papers by Trendafilov and Unkel, in particular Trendafilov \& Unkel (2011) and Unkel \& Trendafilov (2010a\&b), also dealing with the case $p>n$  but proposing alternative models and estimation methods. Trendafilov \& Unkel (2011) appear skeptical to the results of  Robertson \& Symons (2007), and proclaim that when $p>n$  the model assumption of the latter, that $\Psi^2$ is positive definite, is inconsistent with their own model for data. That is certainly right, and we argue below (Sec. 6) that the model for data used by Trendafilov \& Unkel is artificial and unrealistic.  

Our main aim, however, is to show that the fitting of models of type \refeq{Modelforx} in the case of large $p$ is not problematic, and that in any case there is no need to assume normality. We will first derive some basic distribution-free properties of model \refeq{Modelforx}. These are expressed in a normalization of the $x$-components  by $\Psi$,  shown to be suitable for our purpose. It will turn out without difficulties that these properties lead to estimating equations that  are the same as the well-known likelihood equations for $n>p$, thus yielding distribution-free support to the normality-based MLE. 

Another well-known technique for dimension reduction is  principal components analysis (PCA).
PCA aims at describing as much as possible of $\Sigma_{xx}$ by a number of principal components (PCs, linear forms in $x$). There is no model behind PCA, but sometimes the PCs are regarded as representing latent variables in a different, less well-defined way. PCA techniques also have a role in factor analysis. Due to its scale-dependence, the choice of scaling is important.

In the very special case when the error $e$ vanishes, i.e.\ $\Psi^2=0$ in \refeq{Sigmaxx},  $\Lambda \Lambda^T$ can be determined by a PCA on $\Sigma_{xx}$, or estimated by a PCA on the sample covariance matrix $S_{xx}$ (or an SVD on the $x$-data matrix itself). Similarly, if $\Psi^2$ were not zero but regarded as known, we could subtract it from $\Sigma_{xx}$ or $S_{xx}$ and in this way open for use of PCA. This was the basis for the early Principal Factor Analysis method of fitting the FA model: Use some initial $\Psi^2$ to subtract from $S_{xx}$, find PCs yielding an estimate of $\Lambda\Lambda^T$, use this to calculate a new $\Psi^2$, etc. Such methods were found inefficient and unstable, however. In particular they were not scale invariant, in contrast to Gaussian ML  (see Bartholomew \& Knott, 1999, Sec.\ 3.17). From the time when ML methods became computationally feasible and attractive (J\"{o}reskog, 1967, Lawley, 1967), ML estimation has widely replaced the principal factor analysis method. 

In the present paper a new distribution-free method for FA model fitting is proposed, that utilizes principal components of a naturally \emph{rescaled} instead of  \emph{reduced} sample covariance matrix. To our surprise we have not seen this approach in the literature. The methodology has the following properties:
\begin{itemize}
\item It yields the same equations as Gaussian ML--FA for $p<n$, and therefore supports the use of these estimation equations even when the Gaussian distribution is questionable;
\item It is scale invariant in the sense mentioned above;
\item without problems, it allows more variables than observations ($p>n$);
\item It yields estimated or predicted factor scores of high precision when $p$ is large.
\end{itemize}

The basic model properties to be derived in the next section will naturally lead to estimating equations for distribution-free parameter estimation. Different iterative methods to solve these equations are discussed in Section 3. Use of singular value decompositions (SVD) will not only make the computations fast, but also yield some further insight (Sec. 4). The SVD tool is used in Sec.\ 5 to yield expressions for factor scores and residuals.  These are compared in Sec.\ 6 with the model properties of Trendafilov \& Unkel (2011). Finally, in Sec.\ 7, the recommended iteration method is  successfully tried on gene expression data with $p>>n$.

As mentioned above, we assume we have a sample of multivariate $x$-data $x_i$, $i=1,\ldots,n$, $\dim(x)=p$. We will later assume that the $x$-sample is mean-standardized, so we need only consider the sample covariance matrix $S_{xx}=X^T X/(n-1)$ and the corresponding population covariance matrix $\Sigma_{xx}$. In the next section,  we concentrate on $\Sigma_{xx}$, 
so the sample size $n$ and its relation to the dimension $p$ will not yet be a question.

\section{A canonical distribution-free introduction to the FA model}
For the FA model \refeq{Modelforx}, the population covariance matrix $\Sigma_{xx}$ ($p\times p$) is 
\begin{equation}  \label{eq:Sigmaxx}
\Sigma_{xx} = E(S_{xx}) = \Lambda \Lambda^T + \Psi^2 .
\end{equation}
There is a rotational ambiguity in the loading parameters of this representation. For uniqueness we will use the same well-known and natural constraint as in the standard Gaussian ML approach:
\begin{equation}    \label{eq:diagonal}
\Lambda^T \Psi^{-2} \Lambda \quad \textrm{is diagonal}.
\end{equation}
This demand will be equivalent with an assumption that the $p\times k$ matrix 
$\Psi^{-1}\Lambda$ has orthogonal columns. Our motivation to make this particular choice will be clear below. 

As mentioned in Sec.\ 1, classical Principal Factor Analysis requires an initial or current estimate of $\Psi^2$ to be subtracted from  $S_{xx}$, so that ideally we would get $\Lambda \Lambda^T$. PCA is now used on the resulting reduced covariance matrix $S_{xx} - \Psi^2$. Below we will instead use a rescaled covariance matrix, that will be demonstrated to have much better properties.

Consider rescaling the vector $x$ to $z=\Psi^{-1} x$, neglecting for a moment the fact that $\Psi$ is unknown (later we will update $\Psi$ iteratively). This will make all observation components have the same error variance. The total covariance matrix $\Sigma_{zz}$ for a $z$-vector  is 
\begin{equation}  \label{eq:Sigmazz}
\Sigma_{zz} = \Psi^{-1} \Sigma_{xx} \Psi^{-1} = \Psi^{-1}\Lambda (\Psi^{-1}\Lambda)^T + I_p = \Lambda_z \Lambda_z^T + I_p,
\end{equation}
where  $I_p$ denotes the $p\times p$ identity matrix and $\Lambda_z=\Psi^{-1}\Lambda$ is $p\times k$, cf.\ \refeq{Sigmaxx}. Because of assumption \refeq{diagonal} we know that $\Lambda_z$ has orthogonal columns, and it follows that these columns are eigenvectors of the matrix $\Sigma_{zz}$. In a condensed representation we can write 
\begin{equation}  \label{eq:eigenrelation}
\Sigma_{zz} \Lambda_z = \Lambda_z \Omega_z
\end{equation}
where $\Omega_z$ is a diagonal $k \times k$ matrix with the corresponding eigenvalues as diagonal elements, that is 
\begin{equation} \label{eq:Omega} 
\Omega_z = \Lambda_z^T  \Lambda_z + I_k.
\end{equation}
The sum of these $k$ eigenvalues is 
\begin{equation}  \label{eq:traceOmega}
\tr(\Omega_z) =  \tr(\Lambda_z \Lambda_z^T) + k = k + \tr(\Sigma_{zz} - I_p) = k + \sum_{j=1}^p ((\Sigma_{xx})_{jj}-\psi_j^2)/\psi_j^2.
\end{equation}
If $k$ latent factors are both necessary and sufficient for the model to hold, precisely these $k$ eigenvalues of $\Sigma_{zz}$ will be $>1$. For a complete set of eigenvectors of $\Sigma_{zz}$, we need to supplement $\Lambda_z$ by $p-k$ vectors spanning the orthogonal complement of the space spanned by $\Lambda_z$. They will all have the eigenvalue 1.

Equation \refeq{eigenrelation} does not specify the length of the eigenvectors in $\Lambda_z$. For that reason we also introduce  the corresponding set of normalized eigenvectors $\Phi_z$, 
\[
\Phi_z =  \Lambda_z (\Lambda_z^T  \Lambda_z)^{-1/2},
\]
which is a  $p\times k$ matrix of $k<p$ orthonormal eigenvectors. Thus, $\Phi_z^T \Phi_z = I_k$, the  $k\times k$ identity matrix. The matrix $\Phi_z$ of course satisfies the same relation \refeq{eigenrelation} as $\Lambda_z$:
\begin{equation}  \label{eq:Phi-formula}
\Sigma_{zz} \Phi_z = \Phi_z \Omega_z .
\end{equation}

Thus, if we knew $\Psi$ and  $\Sigma_{xx}$, we could form $\Sigma_{zz}$ and calculate its first (=largest) $k$ eigenvectors $\Phi_z$, with their eigenvalues $\Omega_z$, and solve for the loadings matrix $\Lambda=\Psi \Lambda_z$:
\begin{equation} \label{eq:Lambdaz}
\Lambda_z =  \Phi_z  (\Lambda_z^T  \Lambda_z )^{1/2}= \Phi_z \left( \Omega_z - I_k \right) ^{1/2},
\end{equation}
and
\begin{equation} \label{eq:Lambda}
\Lambda = \Psi \Lambda _z = \Psi \Phi_z \left( \Omega_z - I_k \right) ^{1/2}.
\end{equation}
This tells how we can compute $\Lambda$ as a function of $\Psi$ and $\Sigma_{xx}$. In addition, \refeq{Sigmaxx} yields a trivially simple  formula for the diagonal matrix $\Psi^2$  as a function of $\Lambda$, given $\Sigma_{xx}$:
\begin{equation} \label{eq:diagPsi1}
\mathrm{diag}(\Psi^2) = \mathrm{diag}(\Sigma_{xx} - \Lambda \Lambda^T),
\end{equation}
where diag stands for the diagonal part of the matrices, as a vector.
An equivalent alternative is 
\begin{equation} \label{eq:diagPsi2}
\mathrm{diag}(\Psi^{-1} \Sigma_{xx} \Psi^{-1}) = \mathrm{diag}(\Lambda_z \Lambda_z^T + I_p) .
\end{equation}
Here the left hand side can be obtained by elementwise multiplication of the diagonals of $\Psi^{-2}$ and  $\Sigma_{xx}$, or equivalently as 
$\Psi^{-2}\mathrm{diag}(\Sigma_{xx})$.  \\

\section{Parameter estimation}
For parameter estimation based on data, the formulae above can be used with 
$S_{xx}$ inserted for $\Sigma_{xx}$: This yields an estimating equation for $\Lambda$ as
\begin{equation}  \label{eq:Lambdahat}
\widehat{\Lambda}(\Psi)
= \Psi \, \Phi_{z}(S_{zz}) \, \left(\Omega_{z}(S_{zz}) - I_k \right)^{1/2}.
\end{equation}
Here it is indicated that $\widehat\Lambda$ from \refeq{Lambda} is a function of $\Psi$, and that $\Phi_z$ and $\Omega_{z}$ are obtained from $S_{zz}= \Psi^{-1}S_{xx} \Psi^{-1}$ and not from the theoretical $\Sigma_{zz}$. The other estimating equation is obtained from formula \refeq{diagPsi1} or \refeq{diagPsi2} with $S_{xx}$ for $\Sigma_{xx}$:
\begin{equation} \label{eq:diagPsi}
\mathrm{diag}(\widehat{\Psi}^2) = \mathrm{diag}(S_{xx} - \Lambda \Lambda^T),
\end{equation} 
We thus want a solution of these two estimating equations relating $\Psi$ and $\Lambda$.

When $p>n$, these estimating equations turn out to be identically the same as the Gaussian model likelihood equations. This can be taken either as a robustness argument for the Gaussian ML estimates, or as well as a strong argument for the distribution-free method, at least for large $n$. They are also generally quite intuitive. Formula  \refeq{diagPsi} is an obvious demand, and formula \refeq{Lambda} or \refeq{Lambdahat} is a truncated PCA on  $\Lambda \Lambda^T$ after a suitable, albeit parameter-dependent rescaling. 

There is no explicit solution to the set of equations for $\Lambda$ and $\Psi^2$. Thus we have to use some iterative method, and a partial choice is obvious: Select $\Psi^2$ in some way and use this $\Psi^2$ in a calculation of a corresponding $\Lambda$, to be used to update $\Psi^2$, etc. The  step yielding  $\Lambda$ will be taken as given in most of the sequel. The question remains how to update $\Psi^2$. Unless some care is used, such equations  might yield impossible diagonal elements for $\Psi^2$. We return to this question in the next paragraph.

There are alternative estimation methods to ML proposed in the FA literature. Among unweighted and weighted LS metods, the one denoted $\Delta_2$ in Bartholomew \& Knott (1999) appears to be of particular interest in the present context, since it weights data by $\Psi^{-1}$, thus corresponding to our transformation of data.  For given $\Psi$, the  $\Delta_2$ method yields identically the same estimating equation \refeq{Lambdahat}  for $\Lambda$  as the ML method. To estimate $\Psi^2$ by the $\Delta_2$ method is (quoting Bartholomew \& Knott) a good deal more complicated. The choice of $\Psi$ should be such that the sum of squared differences from 1 of the $p-k$ smallest eigenvalues of $\Psi^{-1} S_{xx} \Psi^{-1}$ is as small as possible, under the constraint that they are all $\ge 1$. This constraint, however, excludes the case of a singular $S_{xx}$ and in particular the case $p>n$, and the method is therefore of little interest here.  

Another type of estimation method are the estimation procedures  in for example Trendafilov \& Unkel (2011), jointly estimating $F$, $\Lambda$ and $\Psi^2$. They are based on a different model with additional constraints, which are not adequate in the present setting. They will be further commented in \refsec{U&T}.

\subsection{Iterative solution of the estimating equation system \refeq{Lambdahat} and \refeq{diagPsi}}
The pair of estimating equations  \refeq{Lambdahat} and \refeq{diagPsi}  leads naturally to an iterative procedure, where we start with a provisional  $\Psi$, calculate $\Lambda$ by \refeq{Lambdahat}, calculate a new $\Psi$ by \refeq{diagPsi}, etc. Such calculations are simplified by use of SVD on the sample of $z$-vectors, see next section. However, some variants are possible when using the equations for $\Psi$.

The simplest version is to use \refeq{diagPsi} to express the new $\Psi^2$, in component form 
\begin{equation}  \label{eq:hatpsi1}
\psi_j^2 =(S_{xx})_{jj} -  \left(\Lambda\Lambda^T \right)_{jj}
\end{equation}
with  the current $\Lambda$, based on the previous $\Psi$, on the right hand side. This procedure has a long history, where it turned out to often converge slowly and sometimes to stop before true convergence was achieved. Even worse, the iteration could sometimes yield one or more negative $\Psi^2$ components, known as Generalized Heywood cases. This might be because the best values had not yet been found, but a contributing reason could be the wrong $k$ or an otherwise  inadequate model.  For these reasons, this iteration procedure for Gaussian ML estimation was abandoned, and replaced by a step of direct likelihood maximization to yield $\Psi$ for given $\Lambda$ (J\"{o}reskog, 1967; Lawley, 1967). Another alternative is to use the EM algorithm (Rubin \& Thayer, 1982).

The equivalent formula  \refeq{diagPsi2} suggests a different iteration procedure than \refeq{hatpsi1}. Calculate the new $\Psi$ by \refeq{diagPsi2}, with the current $\Lambda_z$ on the right hand side. This yields the iteration step in component form given by 
\begin{equation}    \label{eq:hatpsi2}
\psi_j^2 =\frac{ (S_{xx})_{jj} }{1 + \left(\Lambda_z\Lambda_z^T \right)_{jj} }.
\end{equation}
One advantage of this is that it yields a positive $\Psi^2$ whatever is the current $\Lambda_z$. On the other hand, our  experiences indicate that it is a slower algorithm, and we do not recommend it.  

Theoretical investigation of the rate of convergence of these methods is difficult, due to the updating of eigenvectors involved. On the other hand, we have used the updating formula  \refeq{hatpsi1} on data with large $p$ ($p>>n$) without any problems, see  further discussion in \refsec{SVD} and \refsec{genedata}. 

\section{Use of the singular value decomposition (SVD)}  \label{sec:SVD}
Let $X$ be the $n\times p$ matrix of column mean-centered $x$-data, and correspondingly $Z=X \Psi^{-1}$ for a provisional $\Psi$. 
A convenient procedure for carrying out the computations above is to calculate and use the singular value decomposition (SVD) of the matrix $Z$, given $\Psi$: 
\[
Z = U D V^T, 
\]
where $U$ ($n\times p$ if $p<n$) and $V$ ($p\times p$) have orthonormal  columns (the left and right singular vectors), and $D$ is a diagonal $p\times p$ matrix whose diagonal elements, the singular values, are, in decreasing order, the square roots of the eigenvalues of  $Z^T Z=V D^2 V^T$. 
When $p>n$, less than $n$ singular values can be positive (typically $n-1$), and then we let $U$ and  $D$ be $n\times n$, and $V$ be $p\times n$.

The right singular vectors forming $V$ are the orthonormal eigenvectors of $Z^T Z$ (or of the covariance matrix $Z^T Z/(n-1)$). Corresponding to the FA model, we truncate the SVD by using only the first $k$ singular vectors, $U_1$ ($n\times k$) and $V_1$ ($p\times k$), say, corresponding to  $\Phi_z$. That is, we partition $Z$ as
\[
Z = U_1 D_1 V_1^T  + U_2 D_2 V_2^T,
\]
where $U=(U_1,\, U_2)$, etc. Note that it does not affect  $U_1 D_1 V_1^T$ whether $p<n$ or $p>n$, but only the second term, where  $D_2$ is either $(p-k)\times (p-k)$ or $(n-k)\times (n-k)$, respectively.

Since $V_1$ is formed by the normalized eigenvectors of $(n-1) S_{zz}$ with the $k$ highest eigenvalues, and these are given by the diagonal $D_1^2$,  we can identify $V_1= \Phi_z$ and $D_1^2=(n-1)\Omega_z$ from equation \refeq{Lambdahat}. Thus the estimating equation \refeq{Lambdahat} for $\Lambda$ can be expressed in terms of $V_1$ and $D_1$, and for the estimation of $\Lambda$ (given $\Psi$) we will need only $U_1 D_1 V_1^T$. More precisely, $\Lambda= \Psi \Lambda_z$ in combination with 
\begin{equation}    \label{eq:hatLambda_z}
\widehat{\Lambda}_z = \Phi_z \left(\Omega_z - I_k\right)^{1/2} = V_1 \left(\frac{D_1^2}{n-1} - I_k\right)^{1/2}.
\end{equation}
Iteration step \refeq{hatpsi1} for $\Psi^2$ takes the following form  in terms of $V_1$ and $D_1$:
\begin{equation}   \label{eq:DiagPsi2New}
\mathrm{diag}(\Psi_{\mathrm{new}}^2) = \mathrm{diag} \left\{ S_{xx} - \Psi \Lambda_z \Lambda_z^T \Psi \right\} =  \mathrm{diag}\left\{S_{xx} - \Psi V_1 \left(\frac{D_1^2}{n-1} - I_k \right) V_1^T  \Psi \right\} .
\end{equation}
The alternative iteration step \refeq{hatpsi2} takes the form
\[
\psi_j^2 =\frac{ (S_{xx})_{jj} }{1 + \left(V_1(D_1^2/(n-1) - 1)V_1^T \right)_{jj} }.
\]

The right hand side of \refeq{DiagPsi2New} may alternatively be expressed as 
\[
\mathrm{diag} \left\{ \Psi \left(V_1 V_1^T + V_2 {D_2}^2 {V_2}^T /(n-1) \right) \Psi \right\} ,
\]
which shows that it is obtained by replacing the first $k$ singular values or eigenvalues  in $S_{zz}$ by the value 1. 
Consequently, the iteration method cannot possibly yield zero or negative values in $\Psi^2$ in any iteration step (presuming start values are positive). What might possibly go wrong, as indicated by  \refeq{hatLambda_z}, is that  $D_1^2/(n-1) - I_k$ is not positive definite.  In the case $p>n$, however, we give below some more results about $D_1^2$ and $D_2^2$, showing that we need not worry.

Note first that when $\Psi$  and $\Lambda$ satisfy the estimating equations, all the $p$ diagonal elements of $S_{zz}-\widehat{\Lambda}_{z} \widehat{\Lambda}_{z}^T$ are 1, so its eigenvalues sum to $p$. At the same time,  
\begin{equation}   \label{eq:residuals1}
S_{zz}-\widehat{\Lambda}_{z} \widehat{\Lambda}_{z}^T= V \frac{D^2}{n-1} V^T - V_1(\Omega_z - I_k) V_1^T =  V_1  I_k V_1^T + V_2 \frac{D_2^2}{n-1} V_2^T.
\end{equation}
Thus, under the same conditions, 
\begin{equation}    \label{eq:p-k}
\tr\left(D_2^2\right)/(n-1) = p-k.
\end{equation}

If $k$ is not higher than motivated by data, we expect the diagonal matrix $\Omega_z - I_k$ in \refeq{hatLambda_z} to have all its diagonal elements positive. When $p<n$, this can fail, and the estimation process too. When $p>n(>k)$, however, the diagonal elements are necessarily positive, at least in a vicinity of the estimation point. To see this, note first that $D_2^2$ contains less than $n-k$ positive values, but has $\tr(D_2^2/(n-1))=p-k$. Thus, the average value is at least $(p-k)/(n-k)>1$. Since the $k$ diagonal values in $\Omega_z=D_1^2/(n-1)$ are larger than this, by selection, the corresponding elements of $\Omega_z - I_k$ are necessarily positive, which was to be shown. 

In passing, we supplement by an expression for the  average of the $k$ first eigenvalues of $S_{zz}$, cf.\ \refeq{traceOmega}. This average can be written
\[
\tr(\Omega_z)/k = 1 +   (\theta - 1) p / k , 
\]  
where $\theta>1$ is the inverse of the harmonic mean of the $p$ unique factor variance proportions $\widehat{\psi}_j^2/(S_{xx})_{jj}$, 
\[
\theta = \frac{1}{p} \sum_{j=1}^p (S_{xx})_{jj} / \widehat{\psi}_j^2.
\] 
This is seen by subtracting $p-k$ from $\tr(S_{zz})$.  Note the proportionality to the dimension $p$ in the second term of $\tr(\Omega_z)$, showing the benefit of large $p$. Note also that when $k$ is increased, $\theta$ will also increase.

\section{Factor scores and model residuals}   \label{sec:scores}
The SVD approach can be used to obtain relatively directly the most common  estimates or predictions of the  scores  $f_i$, or the whole $n\times k$ scores matrix $F$ with the $f$-vectors as rows. As usual in the context of scores estimation/prediction, we provisionally regard the parameters as known (but they are of course estimated). The Bartlett scores, or weighted least squares scores regressing $X$ on $\Lambda$, are given by 
\begin{equation}    \label{eq:Bartlett1}
\widehat{F} = X \Psi^{-2} \Lambda (\Lambda^T \Psi^{-2} \Lambda)^{-1} = Z \Lambda_z (\Lambda_z^T  \Lambda_z)^{-1}
\end{equation}  
so first we can note that with $Z$ as data, Bartlett scores are standard (i.e.\ equal weights) least squares scores. Continuing from  \refeq{Bartlett1},
\[
\widehat{F} = U D V^T \Phi_z (\Lambda_z^T  \Lambda_z)^{-1/2} = U D V^T V_1 (\Omega_z - I_k)^{-1/2} =  U_1 D_1 (\Omega_z - I_k)^{-1/2},
\]
using the fact that $\Phi_z = V_1$.  This implies that the Bartlett score components are proportional to the SVD vectors $U_1$. More precisely, since $D_1^2= (n-1)\Omega_z$, we achieve the following estimation/prediction formula (two equivalent versions related by \refeq{Omega}):
\begin{equation} \label{eq:Bartlett2}
\widehat{F}  = U_1 \sqrt{n-1} \, \Omega_z^{1/2} (\Omega_z - I_k)^{-1/2} = 
U_1 \sqrt{n-1} \left( I_k + (\Lambda_z^T \Lambda_z)^{-1} \right)^{1/2}.
\end{equation}
To the right of $U_1\sqrt{n-1}$ is a diagonal matrix that scales the $j$th column of $U_1$ by the factor $\sqrt{\omega_j/(\omega_j - 1)}$, $j=1,\ldots,k$.
Thus, this is Bartlett's formula in a disguised but computationally convenient form. Typically, if $p$ is large and $k$ is not too large, all  $k$ $\omega$-values will be large (proportionally to $p$, cf. \refeq{traceOmega}), and then with good approximation $\widehat{F} \approx U_1 \sqrt{n-1}$.  

If we instead  predict the scores $F$ by the linear regression of $F$ on the observed $X$-data (or on $Z$),  the best linear predictor $\widetilde{F}$ is given by the so called regression or Thomson scores
\[
\widetilde{F}  = U_1 \sqrt{n-1} \, \Omega_z^{-1/2} (\Omega_z - I_k)^{1/2}  = U_1 \sqrt{n-1} \left( I_k + (\Lambda_z^T \Lambda_z)^{-1} \right)^{-1/2}.
\]
The difference from \refeq{Bartlett2} is the diagonal matrix factor $\Omega_z^{-1}\,(\Omega_z - I_k)$ (cf.\ Bartholomew \& Knott, 1999, sec.\ 3.24, or Krzanowski \& Marriott, 1995, sec.\ 12.27). Again, if $p$ is large, but not $k$,  $\widetilde{F} \approx U_1 \sqrt{n-1}$. 

For high dimension $p$ but small or moderate sample size $n$ we cannot expect high precision in the estimation of $\Lambda$ or $\Psi$. Estimation/prediction of the scores $f_i$, however, will be more precise with increasing $p$.  More precisely, it can be shown that under mild conditions the variance of the factor estimator/predictor $\widehat{F}$ or $\widetilde{F}$ goes to zero as $p$ increases but $k$ and $n$ are kept constant. To be specific, consider the Bartlett score vector $\widehat{f}$  for an arbitrary observation $i$, $\widehat{f}=(\Lambda_z^T \Lambda_z)^{-1} \Lambda_{z}^T z$. 

First, if the difference between $\widehat{\Lambda}_z$ and $\Lambda_z$ is still neglected, formula \refeq{Bartlett1} yields the well-known result
\begin{equation}   \label{eq:varf1}
Var(\widehat{f}\,|\,f) =  (\Lambda_z^T  \Lambda_z)^{-1}  \Lambda_z^T I_p \Lambda_z  (\Lambda_z^T  \Lambda_z)^{-1} =  (\Lambda_z^T  \Lambda_z)^{-1}.
\end{equation}
Due to \refeq{traceOmega}, we may conclude that this diagonal matrix will have small elements when $p$ is large and $k$ is not too large. 

The argument above is not justified when $n<p$, however. In that case, let us still regard $\widehat{\Lambda}_z$ as given, 
but with $\widehat{\Psi}^2$ differing from the right $\Psi^2$. Formula \refeq{varf1} should then replaced by 
\begin{equation}   \label{eq:varf2}
Var(\widehat{f}\,|\,f) =  (\widehat{\Lambda}_z^T \widehat{\Lambda}_z)^{-1} \widehat{\Lambda}_z^T \Psi^2 \widehat{\Psi}^{-2} \widehat{\Lambda}_z  (\widehat{\Lambda}_z^T  \widehat{\Lambda}_z)^{-1} .
\end{equation}
This will differ from the corresponding element of  \refeq{varf1} by less than a factor  
\begin{equation}
\max_j \psi_j^2 / \widehat{\psi}_j^{2}.
\end{equation}
We do not know the true $\psi_j^2$-values, but if there are no components with quite little estimated noise $\widehat{\psi}_j^{2}$, and provided the elements of $(\widehat{\Lambda}_z^T \widehat{\Lambda}_z)^{-1}$ are quite small, we can feel sure the precision in $\widehat{F}$ is high. 

When the scores matrix $F$ has been estimated/predicted, we can form the matrix of residuals, for example $\widehat{E}_x = X-\widehat{F} \widehat{\Lambda}^T$. In order to make them all comparable on the same scale, we must variance-standardize to $\widehat{E}_z=Z-\widehat{F} \widehat{\Lambda}_{z}^T$. Now note that 
\[
\widehat{F} \widehat{\Lambda}_{z}^T =   U_1 D_1 V_1^T  
\]
so the standardized residuals matrix is
\begin{equation}  \label{eq:residuals2}
\widehat{E}_z = Z-\widehat{F} \widehat{\Lambda}_{z}^T =   U_2 D_2 V_2^T . 
\end{equation}
Thus, the sum over $j=1,\ldots,p$ of the mean squared standardized residuals is $V_2 \{D_2^2/(n-1)\} V_2^T$. This may be compared with the result  \refeq{residuals1}, which tells that the trace of $D_2^2/(n-1)$ is only $p-k$, and not $p$, so the mean squared standardized residuals are ``too small'', and must be normalized by $p-k$ instead of $p$ to have the right average size over $j=1,\ldots,p$. This corresponds to the residual degrees of freedom for unbiased variance estimation in a linear model for $Z$, regarding $\Lambda$ as given and the $k(n-1)$ free elements of $F$ as unknowns.  

\section{Models with nonrandom common factors, when $p>n$}   \label{sec:U&T}

In recent years, methods have been advocated for fitting fixed factor models to data, where also $F$ is regarded as a set of unknown parameters, see the review by Unkel \& Trendafilov (2010b). Several papers by those two authors treat the case $p>n$.  The methods of Unkel \& Trendafilov (2010a) and  Trendafilov \& Unkel (2011) proceed from a least squares method minimizing a loss function based on the Frobenius norm of data matrices. Quite generally, the fixed model requires more restrictions than the random model, for uniqueness, and  when $p>n$.  the authors are led to impose special constraints. 
Let us write $X=F \Lambda^T + \Psi E_z$, so we can let $E_z$  exist also when $\Psi^2$ contains zero variances. The papers referred to above assume the model satisfies the constraints  $E_z^T F = 0$, $F^T F \propto I_k$, and (unless $p>n$) $E_z^T E_z \propto I_p$. When  $p>n$ they find that  $E_z^T E_z = I_p$ cannot be fulfilled, because the rank of $E_z$ can be at most $n$, and conclude that they need to allow at least $p-n$ unique factors to have zero variances, corresponding to a singular $\Psi^2$.  In that situation they weaken the constraint $E_z^T E_z = I_p$  to  the eigenvector relation $E_z^T E_z \Psi= \Psi$.

On the other hand, a  result by Robertson \& Symons (2007) states that the Gaussian model likelihood typically (depending on $k$) has a unique global maximum also when $p>n$, and with a \emph{nonsingular} $\Psi$. Trendafilov \& Unkel (2011) correctly remark that  this result is not consistent with their own model.  That the rank of $E_z$ can be at most $n$ (or $n-1$, considering that data are centered) is trivially true for the sample of data, but not for the underlying statistical models assumed by Robertson \& Symons (2007) and by us in the present paper. Our conclusion is that their constraints are artificial, and that their method only represents a constrained partitioning of data, and that it does not represent the fitting of a reasonable statistical model.

We shed further light on this situation here by comparing with our distribution-free but ML-related approach as far as it leads to the eigenvector relation for $\Lambda_z$ and  the Bartlett scores for estimating the scores matrix $F$, with any given $\Psi$: \\
The constraint $E_z^T F = 0$   is satisfied also for the fitted random model and its Bartlett scores $\widehat{F}$, according to \refsec{scores}. \\
The constraint $F^T F \propto I_k$ is not exactly consistent with Bartlett  scores but with the large $p$ approximation $\widehat{F} \approx U_1 \sqrt{n-1}$. \\
The constraints $E_z^T E_z \propto I_p$ for $n>p$ and  $E_z^T E_z \Psi = \Psi$ for $p>n$  are not consistent with our fitted model. and other features of our fitted model, in particular since it does not allow noise outside the diagonal of $E_z^T E_z$.\\
Nor is the constraint consistent with Bartlett  scores and other features of our model. 

As their first illustration, Trendafilov \& Unkel (2011) use Thurstone's 26-variable box data, consisting of a set of $n=20$ boxes and $p=26>20$ variables for each box, representing various aspects of size. When they fit a model with three factors ($k=3$), they get 13 or 14 zero-valued  $\psi_j^2$-values (depending on algorithm). When we fit our model we clearly get no more than 6 zeros, and they can be explained by the peculiarities of the data set. In fact, there are only three original variables in the data set: length, width and height. All other variables are constructed as functions of them. In a model with three latent factors, the factors turn out to be precisely length, width and height, and that explains three zeros. Three other variables are linear functions of length, width and height, and that explains the remaining three zeros. So for example adding a little computer-generated random measurement noise to the variables makes the zero variances disappear completely. Thus, all their zero unique factor variances are not really due to  $n<p$, but to a combination of their assumed artificial data structure (model) and associated fitting method, and the peculiarities of the data set. An example of more applied relevance is studied in \refsec{genedata}.

\section{A gene expression example, with $p=2000$}  \label{sec:genedata} 
We tried the model and the iteration methods on a microarray data set from Alon et al (1999), with 62 tissue samples (a colon cancer sample from each of 40 individuals and non-cancer samples from 22 of these individuals), and $p=2000$ genes (selected by theses authors from a larger set of genes). The data are available on www.bioconductor.org, from where they were fetched. The data have earlier been used for illustrative purposes by McLachlan et al (2003, 2004). 
The response was taken to be the gene expression on log scale (natural log). Each gene was mean- and variance-standardized, but no other normalization of the data was made.  None of the biological structure imposed by the experiment was used in the model, since our aim was not to draw biological conclusions but only to try our methods for model fitting.

We tested the estimation method on the data of  all tissue samples ($n=62$), but mostly on the data of only non-cancer tissue ($n=22$). The iteration method \refeq{hatpsi2} was found to be slower and generally inferior to the method  \refeq{hatpsi1}.
The experiences from running the iteration method \refeq{hatpsi1} were extremely satisfactory.  The method converged in about 10 iterations for small $k$ and not more than 20 to 30 iterations for larger $k$, somewhat also depending on the choice of starting values for $\Psi^2$.  The time per iteration step seemed to be slowly increasing with $k$, but even with an extremely large $k$, $k=20$ say, iterations did not require more time than a second each, on an ordinary laptop. There was no problem of Heywood type during the iterations.  Even if the minimum of the unique factor variances in $\Psi^2$ naturally decreased with $k$, it was in no case estimated to be zero (we tried $k$-values up to 20 for $n=62$, and $k=12$ for $n=22$). After quite few iterations, $\tr(D_2^2)/(n-1)$ was reasonably close to $p-k$, cf. \refeq{p-k}. The statements about $\tr(D_2^2)/(n-1)$ and about the minimum  of the unique factor variances are illustrated in Figures 1 and 2 below, showing how these quantities rapidly converge as the iteration number increases. Both for a small factor dimension ($k=2$) and a moderate ($k=5$) or  large such dimension ($k=12$) there are no problems at all, but $k=10$ is also included for the little bump it shows in Figure 1. 
Starting values were $\psi_j^2=1/2$ for all $j$.

We have thus found substantial support for the conjecture, that the iteration method works so well not \emph{despite} the large $p$-value, but \emph{due to} the large $p$.

\section{Conclusions}
Summing up, we have come to the following conclusions from the investigations in this paper. 

Distribution-free estimating equations for the parameters of the standard FA model (with random factors), \refeq{Lambdahat} and \refeq{diagPsi},  are easily derived in a set-up where variables are variance-normalized by their specific factor standard deviations ($\Psi$). This theory extends the Gaussian likelihood equations both to distribution-free settings and to the case $p>n$. The estimating equations are conveniently expressed by use of a singular value decomposition (SVD) under the same normalization. 

An iteration scheme that has been much used for MLE computation when $p<n$, but also criticized as unreliable in such cases, is shown to have much stronger properties when $p>n$. The theoretical results are supported empirically in an illustration with $p>>n$, where the method was seen to converge quite rapidly.

Another result for situations of type  $p>>n$ is that even though the model parameters cannot be precisely estimated when $n$ is small, the factor scores can be precisely estimated/predicted when $p$ is large.

\section*{References}
Alon, U. et al. (1999). Broad patterns of gene expression revealed by clustering analysis of tumor and normal colon tissue probed by oligonucleotide arrays. \emph{Proc. Nat. Acad. Sci. USA} {\bf 96}, 6745--6750. \\
Bartholomew, D.J. \& Knott, M. (1999). \emph{Latent variable models and factor analysis}, 2nd edn. Arnold, London \\
J\"{o}reskog, K.G. (1967). Some contributions to maximum likelihood factor analysis. \emph{Psychometrika} {\bf 32}, 443--482.  \\
Krzanowski, W.J. \& Marriott, F.H.C. (1995). \emph{Multivariate analysis, part 2}. Arnold, London \\
Lawley, D.N.  (1967). Some new results in maximum likelihood factor analysis. \emph{Proc. Roy. Soc. Edinburgh A} {\bf 67}, 256--264.\\
McLachlan, G.J., Peel, D. \& Bean, R.W. (2003). Modelling high-dimensional data by mixtures of factor analyzers.  \emph{Comp. Stat. \& Data Analysis}  {\bf 41}, 379--388.\\
McLachlan, G.J., Do, K.-A. \& Ambroise, C. (2004). \emph{Analyzing microarray gene expression data}. Wiley, Hoboken.  \\
Rubin, D.B. \& Thayer, D.T. (1982). EM algorithms for ML factor analysis. \emph{Psychometrika} {\bf 32}, 443--482.  \\
Robertson, D. \& Symons, J. (2007). Maximum likelihood factor analysis with rank-deficient sample covariance matrices. \emph{Journal of Multivariate Analysis} {\bf 98}, 813--828. \\
Trendafilov, N.T. \& Unkel, S. (2011).  Exploratory factor analysis of data matrices with more variables than observations. \emph{J. Comp. Graph. Stat.} {\bf 20}, 874--891.\\
Unkel, S \& Trendafilov, N.T. (2010a). A majorization algorithm for simultaneous parameter estimation in robust exploratory factor analysis. \emph{Comp. Stat. \& Data Analysis}  {\bf 54}, 3348--3358.\\
Unkel, S \& Trendafilov, N.T. (2010b). Simultaneous parameter estimation in exploratory factor analysis: an expository review. \emph{Int. Stat. Rev.}  {\bf 78}, 363--382. \\

\vspace{1cm}
\noindent Addresses: \\
Rolf Sundberg, Mathem. statistics, Stockholm University, Sweden, rolfs@math.su.se; \\
Uwe Feldmann, Medical biometry, University of Saarland, Germany, uf@med-imbei.uni-saarland.de \\
Corresponding author: Rolf Sundberg

\begin{figure}[p]
\centerline {
\includegraphics[width=8cm, height=8cm]{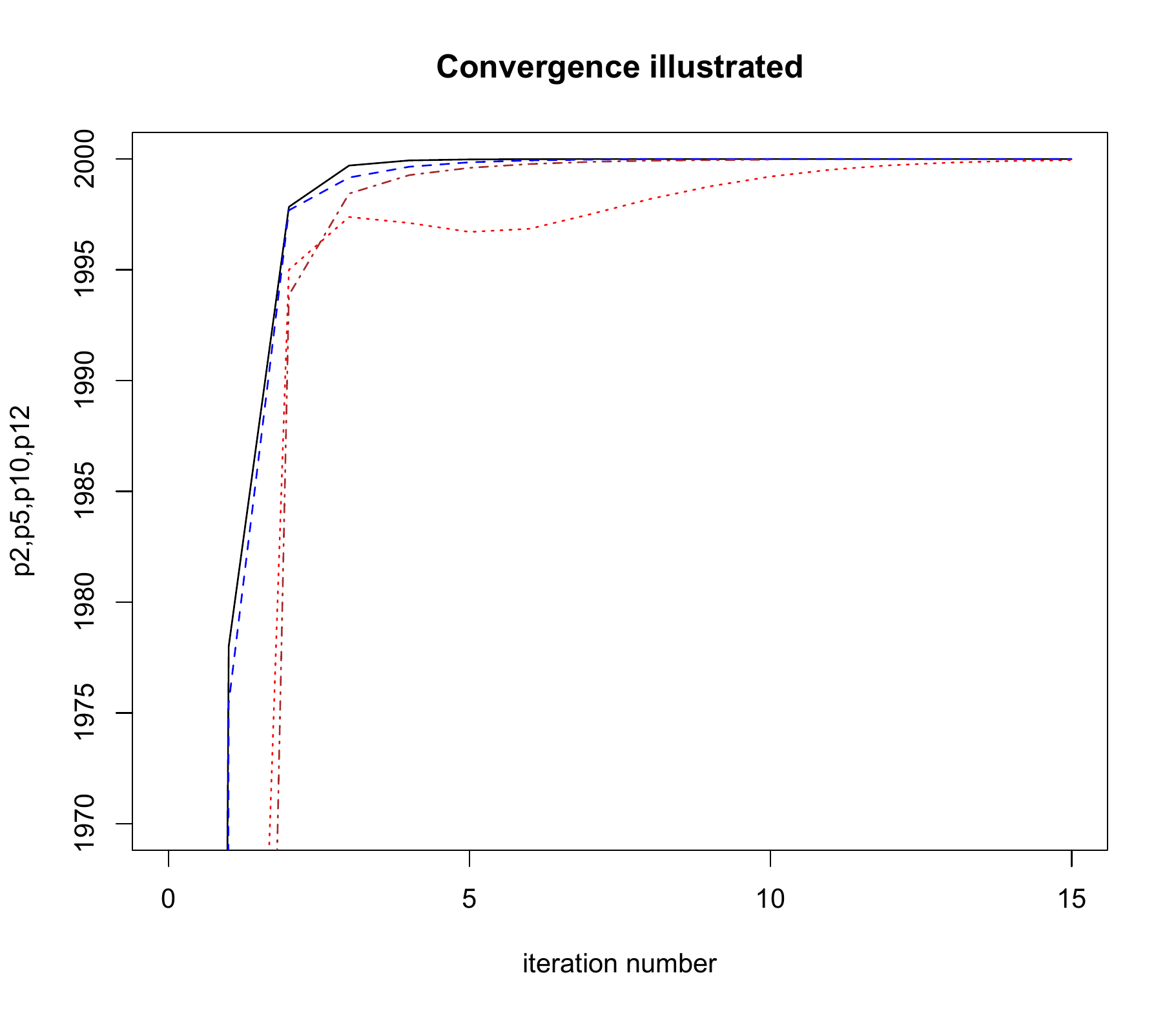}  
}
\caption{Illustrated convergence to $p=2000$ of the sum of unique factor eigenvalues, $+k$, eq.\  \refeq{p-k}.\newline
$k=2$:  ---------  (black)\newline
$k=5$:  - - - - - -  (blue)\newline
$k=10$: $\cdots\cdots$  (red)\newline
$k=12$:  - $\cdot$ - $\cdot$ - (brown)
}
\end{figure}
\begin{figure}[p]
\centerline {
\includegraphics[width=8cm, height=8cm]{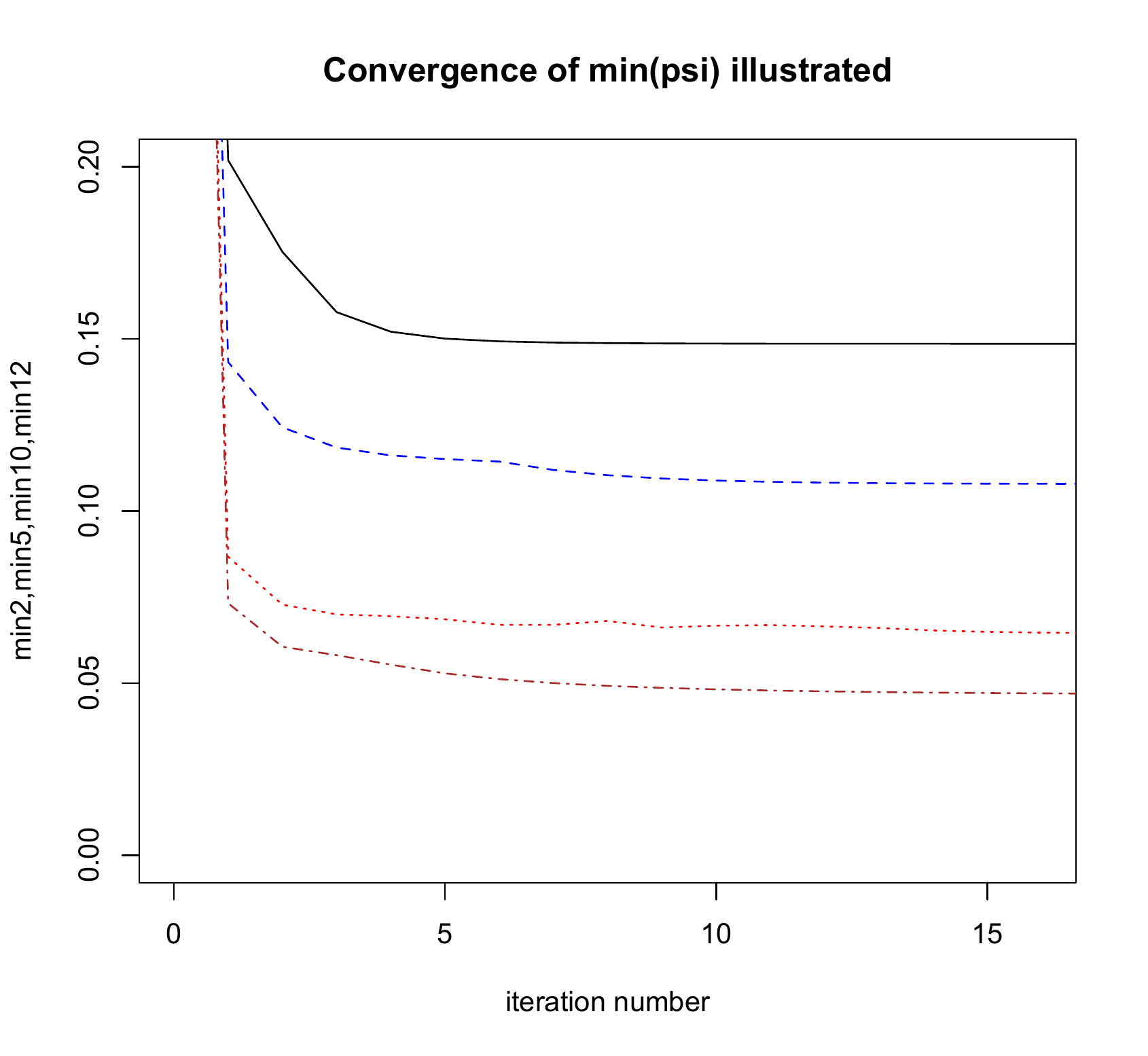}  
}
\caption{Illustrated  convergence                                                                                                                                                                                                                                                                                                                                                                                                                                                                                                                                                                               
of the minimum element of $\Psi$, as the iteration number increases.
$k=2$: ---------  (black)\newline
$k=5$:  - - - - - -  (blue)\newline
$k=10$: $\cdots\cdots$  (red)\newline
$k=12$:  - $\cdot$ - $\cdot$ - (brown)
 }
\end{figure}

\end{document}